\documentclass[12pt]{amsart}

\usepackage{soul}

\usepackage{DLdef1}

\begin{document}

\title[Maximal subalgebras of linear Lie superalgebras]{Maximal
subalgebras of the classical linear Lie superalgebras}

\author{Irina Shchepochkina}
\address{Independent University  of Moscow,
Bolshoj Vlasievsky Per, Dom 11, RU-119 002 Moscow, Russia;
irina@mccme.ru}

{\renewcommand{\thefootnote}{\fnsymbol{footnote}}

\footnotetext{\kern-19pt{\bf AMS Subject classification:}17A70}

\footnotetext{\kern-19pt{\bf  Keywords and phrases:} maximal
subalgebra, Lie algebra, Lie superalgebra}

\begin{abstract} Dynkin's classification of maximal subalgebras of
simple finite dimensional complex Lie algebras is generalized to {\it
linear} Lie superalgebras.  Namely, the maximal non-simple irreducible
subalgebras of $\fgl(p|q)$, $\fq(n)$, $\fsl(p|q)$, $\fosp(m|2n)$,
$\fpe(n)$, and $\fs\fpe(n)$ are classified.
\end{abstract}

\maketitle

\section*{Introduction}

\ssec{Dynkin's result} In 1951, Dynkin published two remarkable
papers somewhat interlaced in their theme: classification of
semi-simple \cite{D1} and maximal \cite{D2} subalgebras of simple
(finite dimensional) Lie algebras.  These classifications are of
interest {\it per se}; they also proved to be useful in the studies
of integrable systems and in representation theory.

A.~L.~Onishchik and D.~Leites\footnote{I am thankful to D. Leites
for help; RFBR grants 99-01-00245 and 01-01-00490a; NFR for partial
financial support and Stockholm University, where most results of
this paper were preprinted in 1988 (no. 32/1988-15), for
hospitality.} asked me to generalize Dynkin's results to \lq\lq
classical" Lie superalgebras. In this paper (partly preprinted in
\cite{Sh1}, \cite{Sh3}) I try to give the answer in a form similar
to that of Dynkin's result.  Let me first remind Dynkin's result. In
what follows the adjective \lq\lq linear" describes a Lie
(super)algebra $\fg$ whose elements are realized as linear operators
in a linear space $V$.

Let $\fg$ be a simple complex linear Lie algebra, e.g.,
$\fg=\fsl(V)$, $\fo(V)$, or $\fsp(V)$, $\fh\subset\fg$ its maximal
subalgebra.  Then only the following 3 cases can occur:

1) the representation of $\fh$ in $V$ is irreducible, $\fh$ is not
simple.  Dynkin's list is as follows (here the cases $\dim V_2=4$ in
the second and fourth lines, as well as $\dim V_2= \dim V_1=2$ in the
third one, are exceptional because $\fo(4)\cong \fsp (2)\oplus \fsp
(2)$):

%\resizebox*{\textwidth}{!}%
{\protect
$$
\renewcommand{\arraystretch}{1.6}
\begin{tabular}{|c|c|c|c|}
\hline
&$\fh$&$\fg$&condition\\
\hline
1&$\fsl (V_1)\oplus \fsl (V_2)$&$ \fsl (V_1\otimes V_2)$&$\dim V_2 \ge \dim V_1 \ge 2$ \\
\hline
2&$\fsp (V_1)\oplus \fo (V_2)$&$\fsp (V_1\otimes V_2)$&
$\dim V_1\ge 2,\;  \dim V_2\geq 3$, $\dim V_2\neq 4$\\
&&&or $\dim V_1=2$ and  $\dim V_2=4$; \\
\hline
3&$\fsp (V_1)\oplus \fsp (V_2)$&$\fo (V_1\otimes V_2)$&$\dim V_2\ge
\dim V_1\ge 2$ \\
&&&except $\dim V_1= \dim V_2= 2$,\\
\hline
4&$\fo(V_1)\oplus \fo (V_2)$&$\fo (V_1\otimes V_2)$&$\dim V_2\ge \dim
V_1\ge 3$\\
&&& and $\dim V_1,\, \dim V_2\neq 4$. \\
\hline
\end{tabular} \eqno{(0.1)}
$$
}

\vskip 0.2 cm

2) $\fh$ is simple and irreducible (i.e., $\fh$ irreducibly acts on
$V$). {\sl  For practically every irreducible representation of a simple
$\fh$ in a linear space $V$ the image is a maximal subalgebra in one
of the three classical simple linear algebras $\fsl(V)$, $\fo(V)$ or
$\fsp(V)$.}  Dynkin proved this and listed the exceptional cases.

3) $\fh$ is reducible.  Then $\fh$ can be described as the
collection of all operators from $\fg$ that preserve a subspace
$W\subset V$. Here $W$ can be arbitrary for $\fg=\fsl$ whereas for
$\fg=\fo(n)$ and $\fsp(2n)$ the bilinear form $\omega$ on $V$
preserved by $\fg$ must either be non-degenerate or identically
vanish on $W$.  Such algebras $\fh$ are certain {\it parabolic}
subalgebras of $\fg$.

Passing to subsuperalgebras we encounter the same cases.  Difficulties in
their superization range widely:

\underline{Superization of case 3)} to Lie superalgebras with Cartan
matrix and linear ones is more or less straightforward, see
\cite{ZZO}.

\underline{Superization of case 2)} requires a more or less explicit
description of finite dimensional irreducible modules over simple
Lie superalgebras.  There are two types of such modules {\it
typical} and {\it atypical} ones, cf.  \cite{K1}, \cite{K2},
\cite{PS}.  Observe that Lemma 4.4 demonstrates that, unlike Lie
algebra case, the images of Lie superalgebras $\fsl(m|n)$,
$\fosp(2|2n)$, $\fpe(n)$, $\fs\fpe(n)$, $\fvect(0|n)$,
$\fsvect(0|n)$, $\fh(0|n)$, $\fh'(0|n)$ in the typical modules are
not maximal, except for the case considered in Lemma 3.4.1.

The case of other algebras and atypical modules constitutes an open
problem.  For a partial result see \cite{J1}.  With a general
character formula, even conjectural, \cite{PS}, \cite{Se} one can now
hope to be able to derive the complete result.

\underline{Superization of case 1)} is what is done in this paper:
the description of the irreducible non-simple maximal
subsuperalgebras of linear complex Lie superalgebras either simple
or \lq\lq classical", i.e., certain algebras closely related to
simple ones.

Regrettably, a precise description of the maximal Lie superalgebras of
type 1) is more involved than Dynkin's description $(0.1)$ above.
Indeed, there are too many exceptional cases occasioned by small
dimensions.  Nevertheless, in Main Theorem, I distinguish four main
types of subalgebras such that any type 1) linear Lie superalgebra is
contained in one of these four types of Lie superalgebras.  Thus, the
subalgebras distinguished in Main Theorem are the main candidates for
the roles of maximal subalgebras.

Observe that two of these four types are similar to Dynkin's types
whereas the other two are of totally different nature, and the picture
is similar to that over fields of prime characteristic.

In \S 1 I describe the main constructions and give a precise
formulation of Main Theorem.  Statements describing when the
subalgebras from Main Theorem are indeed maximal are collected in
Tables 1--3.  \S 2 is devoted to a proof of maximality of
Dynkin-type subalgebras, in \S 3 the other two types are considered.
In \S 4 I prove Main Theorem.

In what follows, $\fii\subplus\fa$ designates a semi-direct sum of
algebras of which $\fii$ is an ideal; same notation is used for
indecomposable modules with a submodule $\fii$.

\ssec{Comparison with the case of prime characteristic} Our results
resemble Ten's result in prime characteristic \cite{T}.  To
formulate it, recall that a subalgebra of a (finite dimensional) Lie
algebra is called {\it regular} if it is invariant with respect to a
maximal torus.

\begin{Theorem} Let $\Kee$ be an algebraically closed
field of characteristic $p>3$.

{\em 1)} Any non-semi-simple maximal subalgebra of $\fsl (n)$ for
$n\not\equiv 0 \mod p$, $\fo (n)$ or $\fsp (2n)$ is regular.

{\em 2a)} Let $\dim V=np^{m}$, $(n, p)=1$, $n>1$.  If $\fm$ is an
irreducible maximal subalgebra in $\fsl(V)$ such that $\fpm =\fm
/\Kee 1_V$ is not semi-simple, then $V=U\otimes O_{m}$, where
$O_{m}=\Kee [x_{1}, \dots , x_{m}]/(x^{p}_{1}, \dots , x^{p}_{m})$
and $\fm =\fgl (U)\otimes O_{m}\subplus \fvect (m)$.

{\em 2b)} If $n=1$, then in addition to the above examples $2a)$ the
algebra $\fm =\fsp(2m)\supplus \fhei (2m|0)$ (where $\fhei (2m|0)$
is an even version of the Heisenberg Lie algebra, see sec.  $2.7$)
is also maximal in $\fsl (V)$.

{\em 3)} Any maximal subalgebra in $\fg _{2}$ is regular except
$\fvect(1)$ for $p=7$ and $\fsl(2)$ for $p>7$.
\end{Theorem}

\ssec{Related results} {\bf Maximal solvable Lie subsuperalgebras}.
Such subalgebras for $\fgl(m|n)$ and $\fsl(m|n)$ are classified in
\cite{Sh2}.  A bizarre series of subalgebras was discovered. Maximal
solvable subalgebras --- Borel subalgebras --- of simple Lie
superalgebras are important in representation theory (e.g., for
construction of Verma modules).  In super setting, the maximal
solvable subalgebras can be larger than what is used for
construction of Verma modules and what Penkov justly suggested to
call Borel subalgebras.  Conjecturally, these larger algebras are
related with atypical representations.

\underline{Superization: case 4)} A nonhomogeneous with respect to
parity subalgebra $\fh$ of the Lie superalgebra $\fg$ is called {\it
Volichenko algebra}.  A list of simple finite dimensional Volichenko
subalgebras in simple Lie superalgebras is obtained under a technical
condition by Serganova \cite{S}.  (For motivations and infinite
dimensional case see \cite{LS}, \cite{KL}.)  Simple Volichenko
algebras are one more, new, type of maximal subalgebras of simple Lie
superalgebras: Volichenko subalgebras are not Lie subsuperalgebras.

\section{Notation, background and main statements}

Before we formulate our result, we need to fix several notations and
constructions.  Formulas of linear algebra are generalized to linear
superalgebra by means of linearity and Sign Rule, such are, for
example, definitions of supercommutator, Lie superalgebra.  There
are, however, notions, e.g., supertrace, which though follow from
Sign Rule, are not obvious direct corollaries, cf. \cite{D}. Some
facts, like existence of (at least) two analogs of the general
linear Lie algebra, or two types of bilinear forms (even and odd)
are even less familiar.

\ssec{Basics} Throughout the paper the ground field is $\Cee$. A
{\it superspace} is a $\Zee/2$-graded linear space $V=V_{\bar
0}\oplus V_{\bar 1}$, where $\Zee/2=\{{\bar 0},{\bar 1}\}$ and
$p(v)=\bar i$ if $v\in V_{\bar i}$.  The {\it superdimension} of $V$
is a pair $N=m|n$, where $m=\dim V_{\bar 0}$, $n=\dim V_{\bar 1}$.
The usual formula $\dim V\otimes W=\dim V\cdot\dim W$ becomes
manifest if we introduce a formal symbol $\eps$ such that $\eps^2=1$
and set $\dim V=\dim V_{\bar 0}+\dim V_{\bar 1}\eps$.

For a superspace $V=V_{\bar 0}\oplus V_{\bar 1}$ denote by $\Pi (V)$
another copy of the same superspace: with the shifted parity, i.e.,
$\Pi(V_{\bar i})= V_{\bar i+\bar 1}$.  The {\it subsuperspace}
$U\subset V$ is a subspace such that $U=U\cap V_{\bar 0}\oplus U\cap
V_{\bar 1}$.  A superspace structure in $V$ induces the superspace
structure in the space $\End (V)$.

The {\it Lie superalgebra} of all linear operators in $V$ is called
the {\it general Lie superalgebra}.  It is denoted by $\fgl(V)$ or
$\fgl(\dim V)$.  Having selected a homogeneous basis of $V$, we can
represent operators by supermatrices; in this paper I only need
supermatrices in the standard format, i.e., when the even basis
vectors of $V$ are collected together and come first.

The space of operators with zero supertrace constitutes the {\it
special linear} Lie subsuperalgebra $\fsl(V)$ also denoted $\fsl(\dim
V)$.

There are, however, at least two super versions of $\fgl(V)$, not one.
Another version is called the {\it queer} Lie superalgebra and is
defined as the one that preserves the complex structure given by an
{\it odd} operator $J$, i.e., is the centralizer $C(J)$ of $J$:
$$
\fq(V)=C(J)=\{X\in\fgl(V)\mid [X, J]=0\},
\text{ where } J^2=-\id.
$$
It is clear that by a change of basis we can reduce $J$ to the form
$J_{2n}=\begin{pmatrix}0&1_n\\ -1_n&0\end{pmatrix}$.  Then in matrix
form we have
$$
\fq(n)=\left \{X\in\fgl(n|n)\mid X=\begin{pmatrix}A&B\\
B&A\end{pmatrix}\right\}.
$$
On $\fq(n)$, the {\it queer trace} is defined: $\qtr:
\begin{pmatrix}A&B\\
B&A\end{pmatrix}\mapsto \tr B$.  Denote by $\fsq(n)$ the Lie
superalgebra of {\it queertraceless} operators.

Observe that the identity representations of $\fg=\fq(V)$ and
$\fsq(V)$ in $V$, though irreducible in super setting, are not
irreducible in the non-graded sense: take linearly independent
vectors $v_1$, \dots , $v_n$ from $V_{\bar 0}$; then
$\Span(v_1+J(v_1), \dots , v_n+J(v_n))$ is a $\fg$-invariant
subspace of $V$ which is not a subsuperspace.

A representation is called {\it irreducible of $G$-type} if it has no
invariant subspace; it is called {\it irreducible of $Q$-type} if it
has no invariant sub{\it super}space, but has an invariant subspace.

\ssec{The action of $\fg_1(V_1)\oplus\fg_2(V_2)$ in $V_1\otimes
V_2$} Given two irreducible (of $G$- or $Q$-type) linear Lie
superalgebras $\fg_i\subset\fgl(V_i)$, we obtain a representation of
the Lie superalgebra $\fg_1\oplus\fg_2$ in the superspace
$V=V_1\otimes V_2$, the {\it tensor product} of the given
representations:
$$
X_1+X_2\mapsto X_1\otimes 1+1\otimes X_2\quad \text{for }X_1\in\fg_1,
\; X_2\in\fg_2.  \eqno{(1.1)}
$$
{\bf G-construction: not both the $\fg_i$ are of $Q$-type}. If both
$\fg_1$ and $\fg_2$ contain the identity operators, the representation
$(1.1)$ has a 1-dimensional kernel.  By $\fg_1\bigodot \fg_2$ we will
mean the image of the direct sum $\fg_1\oplus\fg_2$ under the
representation $(1.1)$.  Observe that $\fg_1\bigodot \fg_2$ is
irreducible in this case.

It is convenient to retain the notation $\fg_1\bigodot \fg_2$ even in
the absence of the kernel, i.e., when $\fg_1\bigodot \fg_2 \cong
\fg_1\oplus\fg_2$.

{\bf Q-construction: both the $\fg_i$ are irreducible of $Q$-type}. In
this case the standard action of $\fg_1\oplus\fg_2$ in $V_1\otimes
V_2$ is reducible.  That is why the construction of $\fg_1\bigodot
\fg_2$ is different, namely, as follows.

Let $U$ be a fixed $(1|1)$-dimensional superspace.  The operators
$J_U= \begin{pmatrix}0 & 1 \\ -1 & 0\end{pmatrix}, \quad I_U=
\begin{pmatrix}0 & i \\ i & 0\end{pmatrix}\in \End(U)$ supercommute
and $J_U^2=I_U^2=-1$.  Let us represent $n_1|n_1$-dimensional
superspace $V_1$ as a tensor product $V_1=(V_1)_{\bar 0}\otimes U$ and
$n_2|n_2$-dimensional superspace $V_2$ as a tensor product
$V_2=U\otimes (V_2)_{\bar 0}$.  Set $J=1\otimes J_U\in\End(V_1)$ and
$I=I_U\otimes 1\in\End(V_2)$.  It is clear that
$$
\renewcommand{\arraystretch}{1.4}
\begin{array}{l}
C(J)=\End(V_1)_{\bar 0}\otimes \Span(1, I_U)\cong \fq(V_1);\\
C(I)= \Span(1, J_U)\otimes \End(V_2)_{\bar 0}\cong \fq(V_2).
\end{array}
$$
Define the {\it $Q$-tensor product} by setting
$$
V_1\otimes^{Q} V_2: =  (V_1)_{\bar 0}\otimes U \otimes (V_2)_{\bar
0}.
$$
The map
$$
X_1+X_2\mapsto X_1\otimes 1_{(V_2)_{\bar 0}}+1_{(V_1)_{\bar 0}}
\otimes X_2\quad \text{for }X_1\in\fg_1,
\; X_2\in\fg_2
$$
determines an irreducible representation of the Lie superalgebra
$\fg_1\oplus\fg_2$ in the space $V_1\otimes^{Q} V_2$.  We will
denote the image of $\fg_1\oplus\fg_2$ in this superspace by
$\fg_1\bigodot \fg_2$.  The usual tensor product $V_1\otimes V_2$
considered as a $\fg_1\oplus \fg_2$-module is a direct sum of two
submodules equivalent to $V_1\otimes_Q V_2$.  (When
$\fg_1=\fg_2=\fg$ and, consequently, $V_1= V_2$, Sergeev denotes the
diagonal action of $\fg$ in $V\otimes^{Q}V$ by $2^{-1}(V\otimes V)$,
cf.  \cite{Ser2}.)

\ssec{Projectivization} If $\fg\subset \fgl(n|n)$ is a Lie
subsuperalgebra containing scalar operators then the {\it
projective} Lie superalgebra of type $\fg$ is $\fpg= \fg/\Cee$.  Lie
superalgebras $\fg_1\bigodot \fg_2$ described in sec.  1.2 are
projective.

Projectivization sometimes leads to new Lie superalgebras, for
example: $\fpgl (n|n)$, $\fpsl (n|n)$, $\fpq (n)$, $\fpsq (n)$;
whereas $\fpgl (p|q)\cong \fsl (p|q)$ if $p\neq q$.

\ssec{Lie superalgebras that preserve bilinear forms} We will often
use a general notation $\faut (\omega)$ for the Lie superalgebra
that preserves the non-degenerate bilinear form $\omega$ in the
superspace $V$, i.e., $\faut (\omega)=$
$$
\{X\in\fgl(V) \mid \omega(Xv_1,v_2)+(-1)^{p(X)p(v_1)}
\omega(v_1,Xv_2)=0  \text{ for any } v_1, v_2\in V\}.
$$
If the form $\omega$ is even and supersymmetric, then the Lie
superalgebra $\faut(\omega)$ is called {\it orthosymplectic} and
denoted $\fosp(V)=\fosp(\dim V)$.  Observe that the passage from $V$
to $\Pi (V)$ sends the supersymmetric forms to superanti-symmetric
ones.  That is why we use the notation $\fosp^{sk}$ for the Lie
superalgebra that preserves the {\it superanti}-symmetric form.  We
have an isomorphism $\fosp(V)\cong\fosp^{sk}(\Pi (V))$, but matrix
representations of elements from $\fosp(V)$ and $\fosp^{sk}(\Pi
(V))$ are different.

If the form $\omega$ is odd, then the Lie superalgebra
$\faut(\omega)$ is called, as A.~Weil suggested, {\it periplectic}
and denoted $\fpe^{sy}(n)$ or $\fpe^{sk}(n)$, in accordance with
symmetry of $\omega$.  The passage from $V$ to $\Pi (V)$ sends the
supersymmetric forms to superanti-symmetric ones and establishes an
isomorphism $\fpe^{sy}(n)\cong\fpe^{sk}(n):=\fpe(n)$.

The {\it special periplectic} superalgebra is
$\fspe(n)=\{X\in\fpe(n)\mid \str X=0\}$.

Observe that the map $\chi_{\lambda}: X\mapsto \lambda\cdot \str X$,
where $\lambda\in\Cee,\; \lambda\ne 0$, determines a nontrivial
character of $\fpe(n)$. Denote by $\fpe_{\lambda}(n)$ the image of
$\fpe(n)$ in the representation $\id \otimes\chi_{\lambda}$.

\ssec{Sergeev Lie superalgebra} A.~Sergeev proved that there is just
one nontrivial central extension of $\fspe(n)$.  It exists only for
$n=4$ and is denoted by $\fas$.  Let us represent an arbitrary
element $X\in\fas$ as a pair $X=x+d\cdot z$, where $x=
\begin{pmatrix}A & B \\ C & -A^t\end{pmatrix}\in\fspe(4)$ ($\tr A=0,
B=B^t, C=-C^t$), $d\in{\Cee}$ and $z$ is the central element.  The
bracket in $\fas$ in the matrix form is
$$
\left[x+d\cdot z,
x' +d'\cdot
z\right]=\left[x,x'
\right]+\tr~CC'\cdot z.
$$

\ssec{Heisenberg Lie superalgebra} Denote by $\fhei(0|m)$ the
Heisenberg Lie superalgebra with $m$ odd generators of creation and
annihilation, i.e., the Lie superalgebra with odd generators
$\xi_1$, \dots , $\xi_n$; $\eta_1$,\dots , $\eta_n$ if $m=2n$ or
$\xi_1$, \dots , $\xi_{n-1}$; $\eta_1$,\dots , $\eta_{n-1}$ and
$\theta$ if $m=2n-1$ and an even generator $z$ satisfying the
relations
$$
\renewcommand{\arraystretch}{1.4}
\begin{array}{l}
    [\xi_i, \eta_j]=\delta_{i, j}\cdot z, \; [\xi_i, \xi_j]=[\eta_i, \eta_j]= [z,
\fhei(0|2n)]=0;\\
{} [\theta, \theta]=z, \; [\theta,\xi_i]=[\theta,\eta_j]=0.\end{array}
$$
Irreducible finite dimensional representation of $\fhei(0|m)$ were
first described in \cite{K}, see also \cite{Ser1}.

Each irreducible finite dimensional representation of $\fhei(0|m)$ is
scalar on $z$ and the only one that sends $z$ to the identity operator
is realized in the superspace $\Lambda (n)=\Lambda(\xi)$ for $m=2n$ or
$\Lambda (n)=\Lambda(\xi,\theta)$ for $m=2n-1$ by the formulas:
$$
z\mapsto 1, \quad \xi_i\mapsto \xi_i, \quad\eta_i\mapsto
\partial_{\xi_i}; \quad\theta\mapsto\theta+\partial_{\theta}.
$$

This representation is irreducible of
$G$-type if $m=2n$ and irreducible of $Q$-type if $m=2n-1$.

The normalizer of $\fhei (0|m)$ in $\fgl(\Lambda (n))$ is $\fg
=\fhei (0|m)\subplus \fo (m)$; it acts in the spinor representation
of $\fo (m)$, or, in terms of differential operators: $\fg=\Span (1,
\xi,
\partial, \xi \partial, \partial\partial, \xi\xi)$ for $m=2n$ and
$\fg=\Span (1, \xi, \partial, \theta+\partial_{\theta}, \xi
(\theta+\partial_{\theta}), \partial(\theta+\partial_{\theta}), \xi
\partial, \partial\partial, \xi\xi)$ for $m=2n-1$.

Let $n=3$.  Then $\fg$ is contained in $\tilde\fg = (\fhei
(0|6)\subplus V)\subplus \fo (6)$ (sum as $\fo (6) $-modules), where
the highest weight of the $\fo (6)\simeq \fsl (4)$-module $V$ is
$(2, 0, 0)$, i.e., $\tilde \fg$ is isomorphic to the nontrivial
central extension $\fas$ of $\fspe (4)$.  (Observe, that $\fsl
(4)$-module $\fhei(0|6)$ is the direct sum of the trivial module and
the exterior square of the dual to the standard 4-dimensional
$\fsl(4)$-module.)

\ssec{Densities} Let $\fvect(0|n)=\fder\Lambda(n)$ be the Lie
superalgebra of vector fields on $(0|n)$-dimensional superspace.
Irreducible representations of $\fvect(0|n)$ are described in
\cite{BL}.  The most important for us will be a one-parameter family
of representations $T^{\lambda}$ of $\fvect(0|n)$ in the superspace
$\Vol^\lambda=\Lambda (\xi)\vvol^\lambda(\xi)$ of
$\lambda$-densities. We define it by the formula
$$
T^{\lambda}(D)(f(\xi)\vvol^\lambda)=(D(f)+\lambda\Div D\cdot f)\vvol^\lambda.
$$
The representations $T^{\lambda}$ are irreducible if $\lambda\ne 0,
1$.  The representation $T^0$ determines the action of the Lie
superalgebra of vector fields $\fvect(0|n)$ in the space of functions;
the constants form an invariant 1-dimensional subspace.  The
representation $T^1$ is the dual representation in the space of volume
forms, $\Vol$, it contains an irreducible subspace of codimension
$\eps^n$ spanned by the volume element $\vvol$.  Therefore,
$T^0\cong(\Pi^n(T^1))^*$.

Define the form $\omega_{1/2}$ on $\sqrt{\Vol}$ by the
formula
$$
\omega_{1/2}(f\sqrt{\vvol}, g\sqrt{\vvol})=\int fg\cdot
\vvol.
$$
It is clear that $T^{1/2}(\fvect(0|n))$ preserves $\omega_{1/2}$.

\ssec{The Poisson superalgebra}
%, a third analog of $\fgl(n)$}
On $\Lambda(m)=\Lambda(\Theta_1,\dots,\Theta_m)$, define a Lie
superalgebra structure by setting (the extra minus is convenient for
calculations with weights)
$$
\{f, g\}_{P.b.}=-(-1)^{p(f)}\sum\limits_{j\leq m}\
\pderf{f}{\Theta_j}\ \pderf{g}{\Theta_j}.
$$
Sometimes it is more convenient to re-denote the $\Theta$'s and set
(here $i^2=-1$ and $[m+1/2]$ is the integer part):
$$
\xi_j=\frac{1}{\sqrt{2}}(\Theta_{j}-i\Theta_{n+j});\;
\eta_j=\frac{1}{ \sqrt{2}}(\Theta_{j}+i\Theta_{n+j})\; \text{ for }\;
j\leq n= [m+1/2]\;\; \theta =\Theta_{2n+1}.
$$
In new indeterminates the Poisson bracket is defined by formula (the
summand with $\theta$ only exists for $m$ odd):
$$
\{f,g\}=-(-1)^{p(f)}\left (\sum_{i\le n}(\frac{\partial f}{\partial
\xi_i} \frac{\partial g}{\partial \eta_i}+ \frac{\partial f}{\partial
\eta_i}\frac{\partial g}{\partial \xi_i}) +\frac{\partial f}{\partial
\theta}\frac{\partial g}{\partial \theta}\right ).
$$
This Lie superalgebra is denoted $\fpo(0|m)$; it is a finite
dimensional analog of the Poisson algebra.  It turns into
the general matrix superalgebra  under
quantization, see \cite{LSh}.

On $\fpo(0|m)$, define the grading to be $\gr(f)=\deg f-2$, where $\deg
f$ is the degree of $f$ as an element of $\Lambda(m)$:
$\deg\Theta_i=1$ for any $i$.  Thus,
$\fpo(0|m)=\mathop{\oplus}\limits_{i=-2}^{m-2}\fpo_i$, where
$\fpo_0\cong \fo(m)$ and where $\mathop{\oplus}\limits_ {i<0}\fpo_i$
is isomorphic to $\fhei(0|m)$.

Since trace and queer trace are quantum versions of the integral,
$\fpo(0|m)$ possesses an ideal $\fspo(0|m)$, the special Poisson
superalgebra, of codimension $\eps^m$,
and a 1-dimensional center, the space of constant functions.

The classical isomorphisms of Lie algebras $\fsl(4)\cong \fo(6)$ and
their modules $\Lambda^2(\id_{\fsl(4)})\cong \id_{\fo(6)}$ show that
$\fas$ described in sec.  1.5 can be embedded into $\fpo(0|6)$.  The
embedding sends the central element $z\in\fas$ into $1\in\fpo(0|6)$.

\ssec{Simplicity} The Lie superalgebras $\fsl(m|n)$ for $m> n\geq
1$, $\fpsl(n|n)$ for $n>1$, $\fpsq(n)$ for $n>2$, $\fosp(m|2n)$ for
$mn\neq 0$ and $\fspe(n)$ for $n>2$ are simple, see \cite{K} (as
well as \cite{Kap}, \cite{FK}, \cite{SNR}).

\ssec{Almost simplicity} We say that a Lie superalgebra $\fg$ is
{\it almost simple} if it can be included (non-strictly) between a
simple Lie superalgebra $\fs$ and the Lie superalgebra $\fder \fs$
of the derivations of the latter: $\fs\subset \fg\subset\fder\fs$.

\ssbegin{Theorem}[Main Theorem] $1\degree$ Let $\fg$ be an
irreducible linear Lie superalgebra which is neither almost simple
nor a central extension of an almost simple Lie superalgebra.

Then $\fg$ is contained in one of the following four major types of
Lie superalgebras:

{\em 1)} $\fgl(V_1)\bigodot \fgl(V_2)$;

{\em 2)} $\fq(V_1)\bigodot \fq(V_2)$;

{\em 3)} $\fgl(V)\otimes \Lambda(n)\subplus \fvect(0|n)$;

{\em 4)} $\fhei(0|2n)\subplus \fo(2n)$.

\noindent $2\degree$ Let, in addition to conditions of $1\degree$,
$\fg$ be a subalgebra of $\fq(V)=C(J)$ for some $J$.  Then $\fg$ is
contained in one of the following Lie superalgebras (numbered as in
$1\degree$):

$1\fq)$ $\fq(V_1)\bigodot \fgl(V_2)$ and $J=J_1\otimes 1$, where
$\fq(V_1)=C(J_1)$;

$3\fq)$ $\fq(V_1)\otimes \Lambda(n)\subplus \fvect(0|n)$ and
$J=J_1\otimes 1$, where $\fq(V_1)=C(J_1)$;

$4\fq)$ $\fhei(0|2n-1)\subplus \fo(2n-1)$.

\noindent $3\degree$ Let, in addition to conditions of $1\degree$,
$\fg$ preserve a non-degenerate homogeneous form $\omega$, either
symmetric or skew-symmetric.  Then $\fg$ is contained in one of the
following Lie superalgebras (numbered as in $1\degree$):

$1\omega)$ $\faut(\omega_1)\bigodot \faut(\omega_2)$ and
$\omega=\omega_1\otimes \omega_2$;

$3\omega)$ $\faut(\omega_1)\otimes \Lambda(n)\subplus
T^{1/2}(\fvect(0|n))$ and $\omega=\omega_1\otimes\omega_{1/2}$.
\end{Theorem}

\ssec{Maximal subalgebras from Main Theorem} Tables 1--2 describe
when subalgebras of the form $\fg_1\bigodot\fg_2$ are maximal in a
linear Lie superalgebra $\fg$. These subalgebras are similar to
those that Dynkin described.

\ssec{Table 1}
In this table we assume that $\fg_i\subset \fgl(V_i)$, and
$\dim V_i=N_i=m_i+ n_i\eps\neq 1$ or $\eps$.
\vskip 0.2 cm

%\resizebox*{\textwidth}{!}%
{\protect$$
\renewcommand{\arraystretch}{1.4}
\begin{array}{|c|c|c|c|}
\hline
\fg _{1} & \fg _{2} & \fg & \text{conditions} \\
\hline
\fgl(N_1) & \fgl(N_2) & \fgl(N_1N_2) & N_i\neq 1+\eps;\quad  m_1\neq n_1
\text{ or } m_2\neq n_2\\
\hline
\fgl(N_1) & \fgl(N_2) & \fsl(N_1N_2) & N_i\neq 1+\eps;\quad m_1= n_1
\text{ and } m_2= n_2\\
\hline
\fsl(N_1) & \fsl(N_2) & \fsl(N_1N_2) & N_i\neq 1+\eps;\quad m_1\neq n_1
\text{ or } m_2\neq n_2\\
\hline
\fq(n_1) & \fq(n_2) & \fsl(n_1n_2|n_1n_2) & n_1n_2> 1\\
\hline
\fq(n_1) & \fgl(m_2+\eps n_2) & \fq(n_1(m_2+ n_2)) & n_1\ge 1;
m_2\neq n_2 \\
\hline
\fq(n_1) & \fgl(m_2+ n_2\eps) & \fsq(n_1(m_2+ n_2)) & n_1\ge 1;
m_2= n_2 > 1 \\\hline
\end{array}
$$
}

\vskip 0.2 cm

The case $N_i=1+\eps$ is exceptional because if we identify the
$(1+\eps)$-dimensional superspace $V$ with $\Lambda(1)$, then
$\fgl(V)\cong \Lambda(1)\subplus \fvect(0|1)$. Therefore,
$$
\fg_1\bigodot \fgl(V)\subset \fg_1\otimes\Lambda(1)\subplus
\fvect(0|1)
$$
for any Lie superalgebra $\fg_1$.  The case $n_1n_2=1$ in the fourth
line is exceptional due to the fact that $\fq(1)\bigodot\fq(1)\cong
\fsl(1|1)$.

\ssec{Table 2}

Table 2 describes maximal subalgebras of the form
$\faut(\omega_1)\oplus \faut(\omega_2)$ of the Lie superalgebra
$\fg$ that preserves a non-degenerate bilinear form
$\omega=\omega_1\otimes \omega_2$.  It is clear that if both forms
$\omega_1$ and $\omega_2$ are supersymmetric or skew then $\omega$
is symmetric while if one of the forms $\omega_1$ or $\omega_2$ is
symmetric and the other one is skew then $\omega$ is skew.  We take
into account isomorphisms $\fosp(V)\cong \fosp^{sk}(\Pi(V))$ and
$\fpe^{sy}(V)\cong \fpe^{sk}(\Pi(V))$, see sec.  1.4, and skip the
types of symmetry just to save space.

In this table $N_i=n_i+2m_i\eps$, and $n_im_i\neq 0$.  The conditions
are occasioned by the fact that the identity representations of
$\fo(2)$ and $\fpe(1)$ are reducible; $\fpe(2)$ is an exception
because
$$
\fpe(2)\cong \fsp(2)\otimes \Lambda(1)\, \subplus\,
T^{1/2}(\fvect(0|1)). \eqno{(1.2)}
$$
If we identify the $2|2$-dimensional superspace of the identity
representation of $\fpe(2)$ with $\Cee^2\otimes \Lambda(1)$, we
obtain an embedding
$$
\faut(\omega_1)\bigodot \fpe(2)\subset
\faut(\omega_1\otimes\omega_2) \otimes\Lambda(1)\subplus
T^{1/2}(\fvect(0|1)),
$$
where $\omega_2$ is the standard form in $\Cee^2$ preserved by
$\fsp(2)$.

\vskip 0.2 cm

%\resizebox*{\textwidth}{!}%
{\protect$$
\renewcommand{\arraystretch}{1.4}
\begin{array}{|c|c|c|c|}
\hline
\fg _{1} & \fg _{2} & \fg & \text{conditions}\\
\hline
\fosp (N_{1}) & \fosp (N_{2}) &
\fosp (N_{1}N_{2})& --- \\
\hline
\fo (n) & \fosp (N_{2}) & \fosp (nN_{2})
& n>2,n\neq  4 \\
\hline
\fsp (2n) & \fosp (N_{2}) &
\fosp (2nN_{2})&n\ge 1 \\
\hline
\fpe (n_{1}) & \fpe (n_{2}) &
\fosp (2n_{1}n_{2}|2n_{1}n_{2})& n_{1},\;  n_{2}>2 \\
\hline
\hline
\fosp (n_{1}|2m_{1}) & \fpe (n_{2}) &
\fpe (n_{1}n_{2}+2m_{1}n_{2})& n_2>2,\; n_{1}\neq 2m_{1} \\
\hline
\fosp (2m|2m) & \fpe (n) &
\fspe (4mn)& n>2\\
\hline
\fosp (n_{1}|2m_{1}) & \fpe_{\lambda} (n_{2}) &
\fpe_{\mu} (n_{1}n_{2}+2m_{1}n_{2}),
& n_2>2, \; n_{1}\neq 2m_{1}\\
& & \text{ for } \mu=\frac{\lambda}{n_1-2m_1}
&  \\
\hline
\fo (n) & \fpe (m) & \fpe (nm) & n, m>2,\quad n\neq 4\\
\hline
\fsp (2n) & \fpe (m) & \fpe (2nm)& m>2,\quad n\ge 1 \\
\hline
\end{array}
$$
}

\ssec{Table 3: Non-Dynkin type of subalgebras} The maximal
subalgebras considered in Tables 1--2 are similar to those
considered by Dynkin.  There are, however, maximal subalgebras of
linear superalgebras of totally different nature.  We represent them
in Table 3: $\fg_1$ is a maximal subalgebra in $\fg$; we set $\dim
V_1=N_1=m_1+ n_1\eps$.  In this case we allow $m_1n_1=0$ but, of
course, exclude $m_1=n_1=0$.  In lines 5--7 we assume that $n_1$ is
even and in lines 8-9 we assume that $m_1=n_1$.  As in Table 2 we do
not mention types of symmetry of bilinear forms.  \vskip 0.2 cm

%\resizebox*{\textwidth}{!}%
{\protect\tiny
$$
\renewcommand{\arraystretch}{1.4}
\begin{array}{|c|c|c|c|}
\hline
&\fg _{1} &  \fg & \text{conditions} \\
\hline 1&\fgl(V_1)\otimes\Lambda(n)\subplus \fvect(0|n) &
\fsl(V_1\otimes\Lambda(n)) & N_1\neq 1+\eps; \text{ either } n\neq 1\\
&&&\text{ or } (n=1 \text{ and } m_1=n_1>1) \\
\hline 2&\fgl(V_1)\otimes\Lambda(1)\subplus \fvect(0|1) &
\fgl(V_1\otimes\Lambda(1)) & m_1\neq n_1;\quad N_1\neq 1 \text{ or }\eps \\
\hline 3&\fgl(V_1)\otimes\Lambda(1)\subplus \Cee\cdot\partial &
\fsl(V_1\otimes\Lambda(1)) & m_1\neq n_1;\quad N_1\neq 1 \text{ or }\eps \\
\hline 4&\fq(V_1)\otimes\Lambda(n)\subplus \fvect(0|n) &
\fsq(V_1\otimes\Lambda(n)) & m_1=n_1\ge 1;\quad n\ge 1 \\
\hline 5&\fosp(V_1)\otimes \Lambda(2k)\subplus T^{1/2}(\fvect(0|2k))
&
\fosp(V_1\otimes \Lambda(2k)) & N_1\neq 1,2;\quad k>0 \\
\hline 6&\fosp(V_1)\otimes \Lambda(2k+1)\subplus
T^{1/2}(\fvect(0|2k+1)) & \fspe(V_1\otimes \Lambda(2k+1)) &
N_1\neq 1,2; \text{ either } k>0\\
&&&\text{ or } (k=0 \text{ and } m_1=n_1>1)\\
\hline 7&\fosp(V_1)\otimes \Lambda(1)\subplus T^{1/2}(\fvect(0|1)) &
\fpe(V_1\otimes \Lambda(1)) & N_1\neq 1, 2;
m_1\ne n_1\\
\hline 8&\fpe(V_1)\otimes \Lambda(2k)\subplus T^{1/2}(\fvect(0|2k))
&
\fspe(V_1\otimes \Lambda(2k)) & m_1=n_1>2;\quad k>0 \\
\hline 9&\fpe(V_1)\otimes \Lambda(2k+1)\subplus
T^{1/2}(\fvect(0|2k+1)) &
\fosp(V_1\otimes \Lambda(2k+1)) & m_1=n_1>2;\quad k\ge 0 \\
\hline
10&\fhei(2n)\subplus \fo(2n) & \fsl(\Lambda(n)) & n\ge 2,\quad  n\ne 3. \\
\hline
11&\fhei(2n-1)\subplus \fo(2n-1) & \fsq(\Lambda(n)) & n> 2. \\
\hline
\end{array}
$$
}

\vskip 0.2 cm

\noindent The exceptional cases:

1) If $\dim V_1=1$ or $\eps$, then $ \fgl(V_1)\otimes
\Lambda(1)\subplus\fvect(0|1)\cong \fgl(1|1)$.

2) If $\dim V_1=1+\eps$, then having identified  $V_1$ with
$\Lambda(1)$, we obtain
$$
\renewcommand{\arraystretch}{1.4}
\begin{array}{l}
    \fgl(V_1)\otimes\Lambda(n)\subplus\fvect(0|n)\cong
(\Lambda(1)\subplus\fvect(1))\otimes\Lambda(n)\subplus\fvect(0|n)=\\
\Lambda(1)\otimes\Lambda(n)\subplus(\fvect(1)\otimes\Lambda(n)
\subplus\fvect(0|n))\subset
\Lambda(n+1)\subplus\fvect(0|n+1).\end{array}
$$

3) The isomorphism $(1.2)$ induces the inclusion
$$
\fpe(2)\otimes\Lambda(n)\subplus T^{1/2}(\fvect(0|n))\subset
\fsp(2)\otimes\Lambda(n+1)\subplus T^{1/2}(\fvect(0|n+1)).
$$

4) $n=3$ in line 10; for motivation see sec. 1.6.

\section{Irreducible maximal subalgebras of the form
$\fg_1\bigodot\fg_2$}

In this section we will prove theorems summarized in Tables 1 and 2.

\ssbegin{Theorem} Let $\dim V_{i}=N_i=m_{i}+ n_{i}\eps$ and $N_i\neq
1$ or $\eps$ or $1+\eps$. Let $V=V_{1}\otimes V_{2}$. Then

{\em 1)} $\fgl (V_{1})\bigodot \fgl (V_{2})$ is a maximal subalgebra
in $\fgl (V)$ if $n_{1}\neq m_{1}$ or $n_{2}\neq m_{2}$, otherwise it
is maximal in $\fsl(V)$.

{\em 2)} The following subalgebras are maximal in $\fsl (V)$:

{\em a)} $\fsl (V_{1})\oplus \fsl (V_{2})$ if $n_{1}\neq m_{1}$ or
$n_{2}\neq m_{2}$;

{\em b)} $\fgl (V_{1})\bigodot \fgl(V_{2})$ if $n_{1}=m_{1}$ and
$n_{2}=m_{2}$. \end{Theorem}

\begin{proof} 1) Let $\fg= \fgl (V_{1})\bigodot \fgl (V_{2})$.
Consider $\fgl(V)$ as $\fg$-module.  Then $\fgl(V)\cong\fgl
(V_{1})\otimes \fgl (V_{2})$ and the bracket in $\fgl(V)$ is defined
via
$$
[A\otimes B, C\otimes D]=(-1)^{p(B)p(C)}[A, C]\otimes
BD+(-1)^{p(A\otimes B)p(C)}CA \otimes [B, D].  \eqno{(2.1)}
$$
Observe that, as $\fgl(V_i)$-module, $\fgl(V_i)$ contains only two
nontrivial submodules: $\Cee\cdot \id$ and $\fsl(V_i)$.  Thus, the
minimal $\fg$-submodule $W$ of $\fgl(V)$ larger than $\fg$ is of the
form
$$
W=\fg + \fsl(V_1)\otimes\fsl(V_2). \eqno{(2.2)}
$$
If $m_1\neq n_1$ and $m_2\neq n_2$, then the sum in (2.2) is direct and
$W=\fgl(V)$, i.e., $\fg$ is maximal in $\fgl(V)$.

Let $\Cee\cdot \id_i\cong \fgl(V_i)/\fsl(V_i)$ be a trivial
$\fgl(V_i)$-module.  If $m_1\neq n_1$ but $m_2= n_2$, then
$\fgl(V)/W\cong \fsl(V_1)\otimes \Cee\cdot \id_2$.  If $m_1= n_1$ and
$m_2= n_2$, then
$$
\fgl(V)/W\cong \fsl(V_1)\otimes \Cee\cdot \id_2+\Cee\cdot \id_1\otimes
\fsl(V_2)+\Cee\cdot \id_1\otimes\Cee\cdot \id_2.
$$
Since the spaces $\fsl(V_i)$ are not closed with respect to the
operator product, formula (2.1) demonstrates that in both cases any
subalgebra strictly containning $\fg$ must contain $ \fsl(V_1)\otimes
\fgl(V_2)+\fgl(V_1)\otimes \fsl(V_2)=\fsl(V)$.  To complete the proof
it suffices to observe that in the first case the subalgebra $\fg$ is
not contained in $\fsl(V)$.
\end{proof}

\ssbegin{Theorem} Let $\dim V_{i}=N_i=m_{i}+ n_{i}\eps$, where
$m_1=n_1\ge 1$ and $N_2\neq 1$, or $\eps$, or $1+\eps$.  Let
$V=V_{1}\otimes V_{2}$.  Then the Lie subsuperalgebra $\fg
=\fq(V_{1})\bigodot \fgl (V_{2})$ is maximal in $\fq (V)$ if
$m_2\neq n_2$; it is maximal in $\fsq(V)$ if $m_2=n_2$.
\end{Theorem}

\begin{proof} The proof of this theorem largely repeats that of the
previous one.  Consider $\fq(V)$ as $\fg$-module.  Then
$\fq(V)\cong\fq(V_{1})\otimes \fgl (V_{2})$ and the bracket in
$\fgl(V)$ is defined via (2.1).

Note that for $m_1>2$ the space $\fq(V_1)$ considered as a
$\fq(V_1)$-module, contains two nontrivial submodules, $\Cee\cdot \id$
and $\fsq(V_1)$.  Therefore, the minimal $\fg$-submodule½ $W\subset
\fq(V)$, is of the form $W=\fg+\fsq(V_1)\otimes\fsl(V_2)$.  Since
$\fsq(V_1)$ is not closed with respect to the operator product and
taking into account formula (2.1), we see that any subalgebra
$\fh\subset\fq(V)$ strictly containing $\fg$ should satisfy
$$
\fh\supset\fsq(V_1)\otimes \fgl(V_2)+\fq(V_1)\otimes \fsl(V_2)=\fsq(V).
\eqno{(2.3)}
$$

For $m_1=2$ the space $\fq(V_1)=\fq(2)$, considered as a
$\fq(2)$-module, contains one more nontrivial submodule, $\Cee\cdot \id
\oplus \fsq(2)_{\bar 1}$.  Therefore, the minimal $\fg$-module
$W\subset \fq(V)$ is in this case of the form $W=\fg+\fsq(V_1)_{\bar
1}\otimes\fsl(V_2)$.  Nevertheless, even in this case formula (2.1)
leads to inclusion (2.3).  To complete the proof for $m_1\ge 2$ it
only remains to observe that for $m_2\neq n_2$ the Lie superalgebra
$\fg$ is not contained in $\fsq(V)$.

The case $m_1=1$ is even simpler: $\fq(V)/\fg\cong \fp\fgl(V_2)$.
\end{proof}

\ssec{Embedding $\fq_1\bigodot \fq_2\subset \fsl(V_1\otimes _Q
V_2)$} First, describe $\fgl(V_i)$ as a $\fq(V_i)$-module.

As in sec. 1.2, introduce a  $(1|1)$-dimensional superspace
$U$ and consider the following basis of $\End (U)$:
$$
1_U, \quad J_U= \begin{pmatrix}0 & 1 \\ -1 & 0\end{pmatrix}, \quad
I_U= \begin{pmatrix}0 & i \\ i & 0\end{pmatrix}, \quad D_U=
\begin{pmatrix}1 & 0 \\ 0 & -1\end{pmatrix}.
$$
Let us realize the $n_1|n_1$-dimensional superspace $V_1$ as a tensor
product $V_1=(V_1)_{\bar 0}\otimes U$ and $n_2|n_2$-dimensional
superspace $V_2$ as a tensor product $V_2=U\otimes (V_2)_{\bar 0}$.
Set $J=1\otimes J_U$ and $I=I_U\otimes 1$.  It is clear that
$$
\renewcommand{\arraystretch}{1.4}
\begin{array}{l}
    C(J)=\End(V_1)_{\bar 0}\otimes \Span(1, I_U)\cong \fq(V_1);\\
C(I)= \Span(1, J_U)\otimes \End(V_2)_{\bar 0}\cong \fq(V_2).
\end{array}
$$
Then $\fgl(V_i)$, as a $\fq(V_i)$-module½ is a direct sum of two
reducible but indecomposable modules:
$$
\fgl(V_1)=\fgl(V_1)_{\bar 0}\otimes\Span(1, I_U)\oplus
\fgl(V_1)_{\bar 0}\otimes\Span(D_U, J_U); \eqno{(2.4)}
$$
$$
\fgl(V_2)=\Span(1, J_U)\otimes\fgl(V_2)_{\bar 0} \oplus \Span(D_U,
I_U)\otimes\fgl(V_2)_{\bar 0}. \eqno{(2.5)}
$$
Observe that the first summands in (2.4) and (2.5) are isomorphic to
$\fq(V_i)$, whereas the second ones, as $\fq(V_i)$-modules, are
isomorphic to $\Pi(\fq(V_i))$.

\sssbegin{Theorem} If $n_1n_2>1$, then
$\fg=\fq(V_1)\bigodot\fq(V_2)$ is maximal in $\fsl(V_1\otimes^{Q}
V_2)$. \end{Theorem}

\begin{proof} Formulas (2.4)--(2.5) show that the minimal
$\fg$-submodule $W\subset \fsl(V_1\otimes^{Q} V_2)$ containing $\fg$
is of the form
$$
W=\fg +\fsl(V_1)_{\bar 0} \otimes \fgl(U)\otimes \fsl(V_2)_{\bar 0}.
$$
When we close $W$ with respect to the bracketing we obtain
$\fsl(V_1\otimes^{Q} V_2)$.
\end{proof}

\ssec{The case with the bilinear form} Let a non-degenerate
homogeneous (with respect to parity) supersymmetric or
superanti-symmetric bilinear form $\omega$ be given in a superspace
$V$.  Consider two objects associated with $\omega$:

1) The Lie superalgebra $\faut (\omega)$ of operators preserving
$\omega$:
$$
\faut (\omega )=\{A\in \fgl (V)\mid \omega (Ax,
y)+(-1)^{p(A)p(x)}\omega (x, Ay)=0\}.
$$

2) The space $\sym (\omega)$ of operators supersymmetric with respect
to $\omega$:
$$
\sym (\omega )=\{A\in \fgl (V)\mid \omega (Ax, y)=
(-1)^{p(A)p(x)}\omega (x, Ay)\}.
$$
It is clear that, as a linear space, $\fgl(V)=\faut(\omega)\oplus \sym
(\omega)$.

Set
$$
\fs\faut (\omega)=\faut (\omega)\cap \fsl (V),\quad\quad \ssym
(\omega) =\sym (\omega)\cap \fsl (V).
$$
Note that if $p(\omega)=\ev$, then $\faut (\omega)= \fs\faut
(\omega)$, and if $p(\omega )=\od$, then $\sym (\omega)=\ssym
(\omega)$.

If $\omega$ is even, it determines a canonical isomorphism $V\cong
V^*$. In this case,
$$
\fg=\faut(\omega)=\begin{cases}\fosp(V)\cong \Lambda^2(V)
\text{ and }
\sym(\omega)=S^2(V)&\text{if
$\omega$ is symmetric}\cr
\fosp^{sk}(V)\cong S^2(V)\text{ and }\sym(\omega)=\Lambda^2(V)&\text{if
$\omega$ is skew}.\end{cases}
$$
In both cases the space
$\sym(\omega)$ contains a 1-dimensional subspace
$\Cee\omega$ corresponding to scalar operators.

If, moreover, $\dim V_{\bar 0}\neq \dim V_{\bar 1}$, then
$\sym(\omega)=\Cee\omega\oplus \ssym(\omega)$, and the $\fg$-module
$\ssym(\omega)$ is irreducible.

If $\dim V_{\bar 0}= \dim V_{\bar 1}$, then
$\Cee\omega\subset\ssym(\omega)$ and $\dim
(\sym(\omega)/\ssym(\omega))=1$ and $\ssym(\omega)/\Cee\omega$ is
irreducible.

If $\omega$ is odd, a canonical isomorphism $V^*\cong
\Pi(V)$ implies $V\otimes V^*\cong V\otimes\Pi(V)\cong\Pi(V\otimes V)$.
In this case,
$$
\fg=\faut(\omega)=\begin{cases}\fpe^{sy}(V)\cong\Pi(\Lambda^2(V))\supset
\fs\faut(\omega)=\fspe^{sy}(V)&\text{if $\omega$ is symmetric}\cr
\fpe^{sk}(V)\cong\Pi(S^2(V))\supset
\fs\faut(\omega)=\fspe^{sk}(V)&\text{if
$\omega$ is skew}\end{cases}
$$
and
$$
\sym(\omega)=\ssym(\omega)=\begin{cases}
\Pi(S^2(V))\supset\Cee\omega&\text{if $\omega$ is symmetric}\cr
\Pi(\Lambda^2(V))\supset\Cee\omega&\text{if
$\omega$ is skew}.\end{cases}
$$
In both cases the space $\sym(\omega)/\Cee\omega$ is an irreducible
$\fg$-module if $\fs\faut(\omega)$ is simple, i.e., if $\dim V=n|n$
and $n>2$.

\sssbegin{Lemma} {\em 1)} $[\sym (\omega ), \sym (\omega )]\subset
\faut (\omega )$, but $[\sym (\omega ), \sym (\omega )]$ is not
contained in $\fs\faut (\omega )$ if $p(\omega)=\bar 1$.

{\em 2)} Set $\{A, B\}= AB+(-1)^{p(A)p(B)}BA$.  Then
$$
\renewcommand{\arraystretch}{1.4}
\begin{array}{l}
    \{\faut (\omega ), \faut (\omega )\},\; \{\sym (\omega ),
\sym (\omega )\}\subset \sym (\omega ); \\
\{\faut (\omega ), \sym (\omega )\}\subset \faut (\omega );\\
\faut (\omega _{1}\otimes \omega _{2})=
\faut (\omega _{1})\otimes \sym (\omega _{2})+\sym
(\omega _{1})\otimes \faut (\omega _{2}); \\
\sym (\omega _{1}\otimes \omega _{2})=
\faut (\omega _{1})\otimes \faut (\omega _{2})+\sym
(\omega _{1})\otimes \sym (\omega _{2}).\end{array}
\eqno{(2.6)}
$$
If $p(\omega)=\bar 0$, then the subspace $\ssym(\omega)$
is not closed with respect to  $\{\cdot, \cdot\}$.

$$
\renewcommand{\arraystretch}{1.4}
\begin{array}{ll}
    3)&[A_1\otimes B_1, A_2\otimes B_2]= \frac
12((-1)^{p(A_2)p(B_1)}[A_1,A_2]\otimes
\{B_1,B_2\}+\\
&(-1)^{p(A_2)(p(A_1+p(B_1))}\{A_2,A_1\}\otimes[B_1,B_2]).\end{array}
\eqno{(2.7)}
$$
\end{Lemma}

Proof: direct calculations. \qed

\ssbegin{Theorem} {\em (Cf. Table 2)} Let Lie superalgebras
$\fs\faut (\omega _{1})$ and $\fs\faut(\omega _{2})$ be simple. Then
Lie subalgebra $\fg =\faut (\omega _{1})\oplus \faut (\omega _{2})$
is maximal in $\faut (\omega _{1}\otimes \omega _{2})$ if either
$p(\omega _{1})+p(\omega _{2}) =\ev$ or if $p(\omega _{i})=\ev$,
$p(\omega _{j})=\od$ and $\dim(V_{i})_{\bar 0}\neq \dim
(V_{i})_{\bar 1}$ for $(i, j)=(1, 2)$.  If $p(\omega _{i})=\ev$,
$p(\omega _{j})=\od$ and $\dim (V_{i})_{\bar 0}=\dim (V_{i})_{\bar
1}$, then $\fg $ is maximal in $\fs\faut (\omega _{1}\otimes \omega
_{2})$.
\end{Theorem}

\begin{proof} Formula (2.6) and the description of $\faut(\omega_i)$
and $\sym(\omega_i)$ as $\faut(\omega_i)$-modules immediately imply
that any subalgebra of $\fh\subset\faut(\omega_1\otimes\omega_2)$
containing $\fg$ must also contain at least one of the submodules
$\fs\faut(\omega_1)\otimes\ssym(\omega_2)$ or
$\ssym(\omega_1)\otimes\fs\faut(\omega_2)$.  But then, by (2.7) we see
that $\fh$ must contain both of these modules, hence,
$$
\fh\supset\fs\faut(\omega_1)\otimes\sym(\omega_2)+
\sym(\omega_1)\otimes\fs\faut(\omega_2). \eqno{(2.8)}
$$
If $p(\omega_i)=\bar 0$ for $i=1,2$, then the rhs of (2.8) coincides
with $\faut(\omega_1\otimes\omega_2)$.  If $p(\omega_i)=\bar 1$ for
$i=1,2$, then by bracketing the elements from distinct summands and
taking into account that
$\{\fs\faut(\omega_i),\sym(\omega_i)\}=\faut(\omega_i)$ we again
obtain that $\fh=\faut(\omega_1\otimes\omega_2)$.

Finally, if $p(\omega_2)=\bar 0$ and $p(\omega_2)=\bar 1$, then the
first summand in (2.8) coincides with
$\faut(\omega_1)\otimes\sym(\omega_2)$.  By bracketing the elements
from distinct summands we see that apart from inclusion (2.8) there is
an inclusion $\fh\supset \ssym(\omega_1)\otimes\faut(\omega_2)$, i.e.,
$\fh\supset \fs\faut(\omega_1\otimes\omega_2)$.

This completes the proof when $\dim (V_{1})_{\bar 0}= \dim
(V_{1})_{\bar 1}$.  For $\dim (V_{1})_{\bar 0}\neq \dim (V_{1})_{\bar
1}$ it suffices to observe that $\fg$ is not contained in $\fsl(V)$.
\end{proof}

\section{Irreducible linear maximal subalgebras of
non-Dynkin's form}

\ssbegin{Theorem} If $m=2n$, where $n\ge 2$, $n\ne 3$, then $\fg
=\fhei (0|m)\subplus \fo(m)$ is a maximal subalgebra in
$\fsl(\Lambda (n))$.

If $n=3$, then $\fg \subset\fas\subset\fsl (\Lambda (3))$, and $\fas$
is maximal in $\fsl (\Lambda (3))$.

If $m=2n-1,\; n>1$, then $\fg$ is a maximal subalgebra in
$\fsq(\Lambda (n))$.
\end{Theorem}

\begin{proof}
Let $\fhei=\fhei(0|m)$.  Consider the image of the universal
enveloping algebra
$U(\fhei)$ in the irreducible representation from sec.  1.6.  Due to
irreducibility this image, considered as a Lie superalgebra, coincides
with $\fgl(\Lambda (n))$ for $m=2n$ and with $\fq(\Lambda (n))$ for
$m=2n-1$.

Consider $U(\fhei)_L$ as a filtered Lie superalgebra with respect to
the filtration induced by the natural filtration of the enveloping
algebra, and consider the associated graded one, $\gr(U(\fhei)_L)$.
As is known (\cite{LSh}), $\gr(U(\fhei)_L)$ is isomorphic to the
Poisson superalgebra with its standard grading described in sec.
1.8, $\fpo(0|m)=\oplus_{i=-2}^{m-2}\fpo_i$.  The graded image
$\gr(\fg)$ of the Lie superalgebra $\fg$ coincides with the
non-positive part of $\fpo$:
$$
\gr(\fg)=\oplus_{-2\leq i\leq 0}^{0}\fpo_i.
$$

Since $\fpo(0|m)$ is transitive, $\fpo_1$ is (for $m\ne 6$) an
irreducible $\fo(m)\cong \fpo_0$-module and generates
$\mathop{\oplus}\limits_{i=1}^{m-3}\fpo_i$, the positive part of the
special Poisson subalgebra, $\fspo(0|m)$, we see that $\gr(\fg)$
is maximal in $\fspo(0|m)$ for $m\ne 6$.

To complete the proof, it suffices to observe that
$\fspo(0|2n)=\gr(\fsl(\Lambda(n)))$ and $\fg\subset \fsl(\Lambda(n))$,
whereas $\fspo(0|2n-1)=\gr(\fsq(\Lambda(n)))$ and $\fg\subset
\fsq(\Lambda(n))$.
\end{proof}

\ssec{$V=V_{1}\otimes \Lambda (n)$} Let $\dim V_{1}=(m_1, n_1)$,
$\fn= \fgl (V_{1})\otimes \Lambda (n)$ and $\fg$ the semidirect sum
of the ideal $\fn$ and the subalgebra $\fvect (0|n)$ with the
natural action on the ideal.  The Lie superalgebra $\fg$ has a
natural faithful representation $\rho$ in $V=V_{1}\otimes \Lambda
(n)$ defined by the formulas for any $A\otimes \varphi \in \fn$,
$D\in \fvect (0|n)$, and $v\otimes \psi\in V$ we have
$$
\begin{matrix}
\rho (A\otimes \varphi )(v\otimes \psi )&= &(-1)^{p(\varphi)
p(v)}Av\otimes \varphi
\psi,\\
\rho (D)(v\otimes \psi )&=&(-1)^{p(D)p(v)}v\otimes D(\psi).\end{matrix}
$$
In the sequel, we will always identify elements of $\fg$
with their images under $\rho$. Therefore, we will consider $\fg$
embedded in $\fgl(V)$ which coincides, as a linear space, with $\End(V)\cong
\End (V_{1})\otimes \End(\Lambda (n))$.

\ssbegin{Theorem} {\em 1)} The Lie superalgebra $\fg =\fgl
(V_{1})\otimes \Lambda (n)\subplus \fvect (0|n)$ is a maximal Lie
subsuperalgebra of $\fsl (V)$, where $V=V_{1}\otimes \Lambda (n)$,
in all cases except

{\em a)} $\dim V_{1}=1+\eps$ or

{\em b)} $n=1$ and $\dim V_{1}=m_1+n_1\eps$
with $m_1\neq n_1$.

{\em 2)} If $n=1$ and $\dim V_{1}=m_1+n_1\eps$, where $m_1\neq n_1$
and $m_1+n_1>1$, then $\fg $ is maximal Lie subsuperalgebra of
$\fgl(V)$. In this case the Lie superalgebra $\fsg:=\fgl
(V_{1})\otimes \Lambda (1)\subplus \Span(\partial) $ is maximal in
$\fsl(V)$.
\end{Theorem}

Let us first prove a particular case of Theorem.

\sssbegin{Lemma} If $n>1$, then $\fg=\Lambda(n)\subplus\fvect(0|n)$
is maximal in $\fsl (\Lambda (n))$.
\end{Lemma}

\begin{proof} We will make use of the arguments used in the proof of
Theorem 3.1.  Observe that $\fg$ contains a $G$-type irreducible
subalgebra $\fhei=\fhei(0|2n)=\Span(1, \xi, \partial)$.  This means
that the image of $U(\fhei)$ in $\End(\Lambda(n))$ coincides with
$\End(\Lambda(n))$, and the graded Lie superalgebra associated with
$U(\fhei)_L$, is isomorphic to $\fpo(0|2n)$.  Clearly, subalgebra
$\fsl(\Lambda(n))$ corresponds to $\fspo(0|2n)$.

Let us realize the elements of $\fpo(0|2n)$ by means of generating
functions in $\xi, \eta$:
$$
\xi_i\longleftrightarrow\xi_i,\quad
\partial_i\longleftrightarrow \eta_i.
$$
The image $\gr(\fg)=\mathop{\oplus}\limits_{i\ge -2}\fg_i$ of $\fg$ in
$\fpo(0|2n)$ is the linear space of functions of degree $\leq 1$ in
$\eta$'s.  In particular, $\fg_0$ consists of the elements of
$\fpo(0|2n)_0\cong \fo(2n)$ that preserve the space
$\Span(\xi_1,\dots,\xi_n)$.  By Dynkin's theorem $\fg_0$ is a
maximal subalgebra in $\fo(2n)$.

On the other hand, it is clear that
$$
\fg_i=\{X\in\fg_{i-1}\mid [X,\fg_{-1}]\subset \fg_{i-1}\}\; \text{ for
all $i\ge 1$},
$$
i.e., $\fg$ is a maximal subalgebra in $\fpo(0|2n)$ with a given
non-positive part.  Therefore, any subalgebra $\fh\subset
\fpo(0|2n)$ containing $\fg$ must contain the whole of $\fpo_0$.
Since for $n\ne 3$ the $\fo(2n)$-action in $\fpo_1$ is irreducible,
$\fh$ must contain the whole component $\fpo_1$, hence, the whole of
$\fspo$.

If $n=3$, then $\fg_1$ does not lie in any of the irreducible
$\fo(6)$-submodules of $\fpo_1$.  Therefore, again,
$\fh\supset\fspo$.

In both cases we see that $\fg$ is maximal in $\fsl(\Lambda(n))$.
\end{proof}

\ssec{Proof of Theorem 3.2.1} Let $\dim V_1\ne 1$ or $\eps$ or
$1+\eps$ and let $\fh\subset \fgl(V)$ be a subalgebra containing
$\fg$.  Clearly, to prove both headings of Theorem, it suffices to
show that $\fh\supset \fsl(V)$.

Observe that $\fg$ contains two subalgebras, $\fg_1=\fgl(V_1)\otimes
1$ and $\fg_2=1\otimes\Lambda(n)\subplus\fvect(0|n)$.  Since
$\End(V_1)$ only contains two nontrivial $\fgl(V_1)$-submodules,
$\Cee\cdot \id$ and $\fsl(V_1)$, we deduce that $\fh$ must contain a
subspace of one of the two types: either $1\otimes W$ or
$\fsl(V_1)\otimes W$ for a $\fg_2$-invariant subspace
$W\subset\fgl(\Lambda(n))$.

In the first case, by Lemma 3.2.2, $\fh\supset 1\otimes
\fsl(\Lambda(n))$ and, by Theorem 2.1, $\fh\supset \fsl(V_1\otimes
\Lambda(n))$.

In the second case, let us realize the elements of $\fgl(\Lambda(n))$
by differential operators acting on $\Lambda(n)$.  It is clear
that by bracketing with $\xi_i$ and $\partial_i$ we can reduce any
differential operator to the form $\partial_j$, i.e., $W\supset
\partial_1, \dots, \partial_n$.  Applying formula $(2.1)$ to elements
of the form $A\otimes \partial_j$ and $B\otimes D$, where
$D\in\Lambda(\partial)$ is a differential operator with constant
coefficients, we see that $W\supset \Lambda(\partial)$.  Therefore, by
$\fg_2$-invariance, $W\supset \fsl(\Lambda(n))$ and, therefore,
$\fh\supset \fsl(V_1\otimes \Lambda(n))$.  \qed

\ssbegin{Theorem} The Lie superalgebra $\fg =\fq (V_{1})\otimes
\Lambda (n)\subplus \fvect (0|n)$ is a maximal Lie subsuperalgebra
of $\fsq (V)$, where $V=V_{1}\otimes \Lambda (n)$.
\end{Theorem}

Proof is similar to that of Theorem 3.2.1. \qed

\ssec{Lie superalgebras preserving a non-degenerate bilinear form}
We use notations from sec. 1.4 and 1.7. Observe that for $n\le 2$ we
have isomorphisms:
$$
T^{1/2}(\fvect(0|1))\cong \fpe(1);\quad T^{1/2}(\fvect(0|2))\cong
\fosp(2|2).
$$

\sssbegin{Lemma} For $n>2$ the Lie superalgebra
$\fg=T^{1/2}(\fvect(0|n))$ is maximal in $\fs\faut(\omega_{1/2})$.
\end{Lemma}

\begin{proof} Let us realize the space $\fgl(\Lambda(\xi))$ by
differential operators.  Observe that any function
$\varphi\in\Lambda(\xi)$ is an element from the space
$\sym(\omega_{1/2})$, and any differential operator with constant
coefficients $D\in\Lambda^k(\partial)$ belongs to
$\faut(\omega_{1/2})$ for $k$ odd and to $\sym(\omega_{1/2})$ for $k$
even.  Since $\{A, B\}= 2AB-(-1)^{p(A)p(B)}[A, B]$, it follows that
$$
\renewcommand{\arraystretch}{1.4}
\begin{array}{l}
\faut(\omega_{1/2})=\Span(\varphi, D \mid \varphi\in\Lambda(\xi),\;
D\in\mathop{\oplus}\limits_k \Lambda^{2k+1}(\partial));\\
\sym(\omega_{1/2})=\Span(\varphi, D \mid \varphi\in\Lambda(\xi),\;
D\in\mathop{\oplus}\limits_k \Lambda^{2k}(\partial)).\end{array}
$$
Set $V_l= \mathop{\oplus}\limits_{k\le l,\; k\equiv l\pmod
2}\{\Lambda(\xi),\;\Lambda^{k}(\partial)\}$.  Clearly, each $V_l$ is
$\fg$-invariant, and as follows from \cite{BL}, the $\fg$-action in
$V_l/V_{l-2}$ is irreducible for $l\ne 0, n$ (we assume that
$V_{-2}=V_{-1}=0$).

The explicit formula for the bracket of
$A=\{\xi_1\xi_2,\partial_1\}\in\fg$ and
$B=\{\xi_1\xi_2,\;\partial_1\partial_2\partial_{i_3}\dots
\partial_{i_k}\}\in V_k$ is
$$
[A,B]=-\xi_2\partial_{i_3}\dots
\partial_{i_k}\in V_{k-2};
$$
it shows that the representation of $\fg$ in $V_k$ is not completely
reducible and the minimal $\fg$-invariant subspace $W$ such that
$\fg\subset W\subset\faut(\omega_{1/2})$ is $V_3$ for $n>3$ and
$\fs\faut(\omega_{1/2})\cong \fspe(4)$ for $n=3$.

Finally, the formula
$$
[\partial_1\partial_2\partial_3,\;
\xi_1\partial_{i_1}\partial_{i_2}\dots
\partial_{i_k}]=\partial_2\partial_3\partial_{i_1}\dots\partial_{i_k}
\in V_{k+2}
$$
shows that $V_3$ generates the whole Lie superalgebra
$\fs\faut(\omega_{1/2})$.
\end{proof}

A corollary from the proof of Lemma 3.4.1 is the following statement:
{\sl the minimal $\fg$-invariant subspace $W$ such that
$\Lambda(n)\subset W\subset\sym(\omega_{1/2})$ is $W=V_2$}.

\ssbegin{Theorem} [Cf.  Table 3] Let $\omega_1$ be a non-degenerate
supersymmetric or skew bilinear form in $V_1$ of dimension $N_1=m_1+
n_1\eps$, where $N_1\ne 1,\; 2$ and $n_1$ is even if
$p(\omega_1)=\ev$, and where $m_1=n_1>2$ if $p(\omega_1)=\od$.

Then $\fg=\faut(\omega_1)\otimes\Lambda(n)\subplus
T^{1/2}(\fvect(0|n))$ is a maximal subalgebra in
$\fs\faut(\omega_1\otimes\omega_{1/2})$ except for the case when
$n=1$, $p(\omega_1)=\ev$ and $m_1\ne n_1$; in this case $\fg$ is
maximal in $\faut(\omega_1\otimes\omega_{1/2})$.
\end{Theorem}

Lemma 3.4.1 and its Corollary make it possible, essentially, to combine
proofs of Theorem 3.2.1 and 2.4.2.

\section{Proof of Main Theorem}

\ssec{Superization of a Proposition by Dixmier} In this section
$\fg$ is an irreducible linear Lie superalgebra, $\rho$ its standard
representation in a finite dimensional superspace $V$ (in particular
$\rho$ is faithful).  Let $\fii$ be an ideal of $\fg$.  We will
assume that $\dim V>1$ (because the case $\dim V=1$ is trivial). Our
proof of Main Theorem is largely based on the following
constructions and statements (Theorems 4.1.1 and 4.1.2), a
superization of statements well-known to the reader from Dixmier's
book \cite{Di} (Proposition 5.5.1).

Let $\tau$ be an irreducible subrepresentation of the restriction
$\rho|_{\fii}$ in the subspace $U\subset V$ and let $V_1\subset V$ be
the sum of all $\fii$-submodules of $\rho|_{\fii}$, isomorphic to
either $\tau$ or $\Pi(\tau)$.

The {\it stabilizer} $\fst(\tau)$ of $\tau$ is the set of $Y\in\fg$
such that there exists an $A\in\End(U)$ for which
$\tau([Y,X])=[A,\tau(X)]$ for any $X\in \fii$.

\sssbegin{Theorem} {\em 1)}  $\fh=\fst(\tau)\supset\fg_{\ev}$;

{\em 2)} $V_1$ is  $\fh$-invariant;

{\em 3)}If $\fh\neq \fg$, then $\rho\equiv \Ind^{\fg }_{\fh }\sigma$,
where $\sigma$ is a representation of $\fh$ in $V_1$.
\end{Theorem}

\sssbegin{Theorem} Let $\fg$ be an irreducible linear Lie
superalgebra, $\rho$ its standard representation (in particular,
$\rho$ is faithful). Let $\fii $ be a nontrivial ideal of $\fg$.
Then $3$ cases are possible:

{\em A)} $\rho |_{\fii}$ is irreducible;

{\em B)} $\rho |_{\fii}$ is a multiple of an irreducible $\fii$-module
$\tau$ and the multiplicity of $\tau $ is $>1$;

{\em C)} there exists a proper subalgebra $\fh\subset\fg$ such that
$\fii\subset \fh$ and $\rho\equiv \ind^{\fg}_{\fh }\sigma$ for an
irreducible $\fh$-module $\sigma$.
\end{Theorem}

Proof largely follows the lines of \cite{Di} with one novel case:
irreducible modules of $Q$-type might occur.  Fortunately, their
treatment is rather straightforward and I would rather save paper by
omitting the verification.

We will say that the {\it representation $\rho$ is of type $A$, $B$,
or $C$ with respect to the ideal} $\fii$ if the corresponding case
holds.  Lemmas 4.5, 4.2 and 4.4  deal with types A, B and C,
respectively.

\ssbegin{Lemma} {\em {\bf B)} 1)} If $\rho$ is of type B with
respect to the ideal $\fii$, then either $\fg\subset
\fgl(V_1)\bigodot\fgl(V_2)$ or $\fg\subset \fq(V_1)\bigodot\fq(V_2)$
for some $V_1$ and $V_2$.

{\em 2)} If, moreover, $\fg\subset\fq(V)$, then $\fg\subset
\fq(V_1)\bigodot\fgl(V_2)$ for some $V_1$, $V_2$;

{\em 3)} If $\omega$ is a non-degenerate $2$-form on $V$ and
$\fg\subset \faut(\omega)$, then $\fg\subset \faut(\omega_1)\bigodot
\faut(\omega_2)$, where $\omega=\omega_1\otimes \omega_2$
\end{Lemma}

Proof of Lemma 4.2 follows, except statements involving $\fq$, the
standard scheme of the proof of a similar statement for Lie algebras,
and the exception is also easy to consider, so I skip it.

\ssec{Remark}  If $\rho$ is of type B with respect to the
(nontrivial) ideal $\fii$ and $\dim \tau=1$ or $\eps$ then due to
the faithfulness of $\rho$ the ideal $\fii$ should be a
1-dimensional center of $\fg$.

\ssbegin{Lemma} {\em {\bf C)} 1)} If $\rho$ is of type C with
respect to the ideal $\fii$, then $\fg\subset \fgl(V_1)\otimes
\Lambda (n)\subplus\fvect(0|n)$ for some $V_1$ and $n$ and
$V=V_1\otimes\Lambda(n)$.

{\em 2)} If, moreover, $\fg\subset\fq(n)$, then  $\fg\subset
\fq(V_1)\otimes\Lambda (n)\subplus\fvect(0|n)$

{\em 3)} If $\rho$ is of type C with respect to the ideal $\fii$ and
$\fg\subset \faut(\omega)$ for a non-degenerate form $\omega$, then
$\fg\subset\faut(\omega_1)\otimes \Lambda(n)\subplus
T^{1/2}(\fvect(0|n))$, where $\omega=\omega_1\otimes \omega_{1/2}$
and the representation $T^{1/2}$ of $\fvect(0|n)$ in the superspace
$\sqrt{\Vol}$ of half-densities is defined in sec. $1.7$.
\end{Lemma}

\begin{proof}
By definition of type C, $\rho=\ind_\fh^\fg\sigma$. Let
$$
\fh_1=\{X\in\Ker\sigma\mid[X,\fg]\subset \fh\}.
$$
Clearly, $\fh_1$ is a subalgebra in $\fg$ and ideal in $\fh$.
For each integer $i>1$ define inductively a subalgebra
$$
\fh_i=\{X\in\fh_{i-1}\mid[X,\fg]\subset \fh_{i-1}\}.
$$
We obtain a decreasing filtration in $\fg$:
$$
\fh_{-1}=\fg\supset\fh_0=\fh\supset\fh_1\supset\dots
$$
Since $\dim\fg<\infty$, the filtration stabilizes, i.e.,
$\fh_k=\fh_{k+1}=\dots$ for some $k$.  This means that $\fh_k$ is an
ideal in $\fg$ lying in the kernel of $\sigma$.  But then
$\fh_k\subset\Ker \rho$ and, therefore, $\fh_k=0$ because $\rho$ is
faithful.

Consider the associated graded Lie
superalgebra
$$
\gr(\fg)=\mathop{\oplus}\limits_{i=-1}^{k-1}\fg_i, \text{ where }
\fg_i=\fh_i/\fh_{i+1}.
$$
Observe that we have two homomorphisms of $\fh$:
$$
\sigma: \fh\tto\fgl(V_1) \quad\text{ and }\quad
\ad_{\fg/\fh}:\fh\tto\fgl(\fg/\fh),
$$
and $\fh_1=\Ker\sigma\cap\Ker(\ad_{\fg/\fh})$.  Set $W=\fg/\fh$, and let
$\dim W=0|n$.  We have obtained an embedding
$$
\fg_0\subset\fgl(V_1)\oplus\fgl(W)\cong\fgl(V_1)\oplus W^*\otimes W,
$$
which for every $i>0$ induces an embedding
$$
\fg_i\subset\fgl(V_1)\otimes\Lambda^{i}(W^*)\oplus \Lambda^{i+1}(W^*)
\otimes W,
$$
which add up to an embedding of the whole $\fg$:
$$\gr(\fg)\subset\fgl(V_1)\otimes
\Lambda(W^*)\subplus\fvect(W^*).
$$

It remains to observe that since $\fvect(W^*)$ contains a grading
operator, the embedding exists not only for $\gr(\fg)$, but for $\fg$
itself.  Under the embedding the space $V$ is identified with
$V_1\otimes \Lambda(W)\cong V_1\otimes\Lambda(W^*)^*$.  As we observed
in sec.  1.7, the space $\Lambda(W^*)^*$, as a $\fvect(W^*)$-module½
is isomorphic  $\Pi^n(T^1)$. Since we are interested
not in the representation itself but only in the image of
$\fg$ in $\fgl(V)$, we can replace for convenience
$V_1\otimes \Lambda(W)$ with $\tilde{V_1}\otimes \Lambda(W^*)$ for
some $\tilde{V_1}$. This completes the proof of the first heading of
Lemma.

To prove the second heading, observe that if $\fg\subset
\fq(V)=C(J)$, then subspace $V_1$ is $J$-invariant and,
therefore, $\sigma(\fh)\subset \fq(V_1)=C(J|_{V_1})$.

Finally, if $\fg$ preserves a non-degenerate bilinear form $\omega$
in $V$, consider its restriction $\tilde\omega$ on $V_1$.  Since
$\sigma$ is irreducible, $\tilde\omega$ is either non-degenerate or
vanishes identically.  Denote by $V_1^\perp$ the subspace orthogonal
to $V_1$ with respect to $\omega$.  Clearly, $V_1^\perp$ is
$\fh$-invariant.

If $\tilde\omega$ is non-degenerate, then $V=V_1\oplus V_1^\perp$.
This means that $\fg\subset\faut(\tilde\omega)\otimes 1\oplus
1\otimes \fvect(0|n)$.  Therefore, the $\fg$-action in $V$ is
reducible; contradiction.

If $\tilde\omega\equiv 0$, then by $\fg$-invariance of $\omega$, the
space $V_1^\perp$ must contain the image of $V$ under $\rho(\Ker
\sigma)$.  Hence, the $\fh$-module $V/V_1^\perp\cong V_1^*$ should
be of the form $V_1\otimes \Lambda^n(W)$.  This implies that
$\omega=\omega_1\otimes \omega_{1/2}$ for some $\omega_1$ and
$\fg\subset \faut(\omega_1)\otimes\Lambda(n)\subplus
T^{1/2}(\fvect(0|n))$. \end{proof}

To complete the proof of Main Theorem, it suffices to
consider the case when $\rho$ possesses the following property: for
any ideal $\fii\subset \fg$ either $\rho$ is of type A with respect to
$\fii$ or $\rho\vert _\fii$ is the multiple of a character. Due to
Remark 4.3 the second possibility means that $\fii$ is a 1-dimensional
center of $\fg$. In particular, any nontrivial commutative ideal of $\fg$
coincides with its 1-dimensional center.

Let $\fr$ be the radical of the linear Lie superalgebras $\fg$.
We see that either $\dim \fr\le 1$ or $\fr$ is not commutative.

In the case when $\fr$ is not commutative, consider the
derived series of $\fr$:
$$
\fr\supset\fr_1\supset\dots\supset\fr_k\supset\fr_{k+1}=0,
\text{ where $\fr_{i+1}=[\fr_i, \fr_i]$.}
$$
Clearly, each $\fr_i$ is an ideal in $\fg$ and the last ideal, $\fr_k$
is commutative.  Hence, $\dim \fr_k=1$ and $\fr_k$ is the center of $\fg$.

\ssbegin{Lemma} {\em {\bf A)} 1)} If $\rho$ is of type A with
respect to $\fr_{k-1}$ and $\rho| _{\fr_{k}}$ is scalar, then either
$\fr_{k-1}\cong\fhei(0|2n)$ or $\fr_{k-1}\cong\fhei(0|2n-1)$ and
$V\cong\Lambda(n)$ or $\Pi(\Lambda(n))$ and $\fg\subset
\fhei(0|2n)\subplus\fo(2n)$.

{\em 2)} If additionally $\fg\subset \fq(V)$, then
$\fr_{k-1}\cong\fhei(0|2n-1)$ and
$\fg\subset\fhei(0|2n-1)\subplus\fo(2n-1)$.

{\em 3)} Under assumptions of heading $1)$ $\fg$ does not preserve
any non-degenerate  bilinear form on $V$.
\end{Lemma}

\begin{proof} We will prove headings 1 and 2 simultaneously.

1) Since $\fr_k$ is the center of $\fg$ and $\dim\fr_k=1$
we have $\fr_{k-1}=\fhei (0|m)$ for some $m$;

2) $\rho|_{\fr_{k-1}}$ is irreducible and faithful, so it can be
realized in the superspace of functions, $\Lambda (n)$, or in
$\Pi(\Lambda (n))$, where $n=[\frac {m+1}2]$; observe that
$\rho|_{\fr_{k-1}}$ is irreducible of  $G$-type for $m=2n$ and
it is irreducible of $Q$-type for $m=2n-1$, see sec. 1.6;

3) $\fg$ is contained in the normalizer of $\fhei (0|m)$ in
$\fgl(\Lambda(n))$, i.e., $\fg\subset\fhei (0|m)\subplus \fo(m)$.
\end{proof}

\ssbegin{Lemma} Let $\dim\fr\leq 1$, i.e., $\fg$ is either
semi-simple or a nontrivial central extension of a semi-simple Lie
superalgebra but NOT an almost simple or a central extension of an
almost simple Lie superalgebra.  Then we can always choose an ideal
$\fii\ne \fr$ such that $\rho$ is of type B or C with respect to
$\fii$.
\end{Lemma}

\ssec{Proof of Lemma 4.6 for semi-simple Lie superalgebras}

By definition $\fg$ is {\it semi-simple} if its radical is zero.  By
analogy with description of semi-simple Lie algebras over fields of
prime characteristic, V.~Kac \cite{K} described semi-simple finite
dimensional Lie superalgebras as follows.  Let $\fs_1$, ...  ,
$\fs_k$ be simple Lie superalgebras, let $n_1$, ...  , $n_k$ be
nonnegative integers, $\Lambda(n_j)$ be the supercommutative
Grassmann superalgebra, and $\fs=\oplus\fs_j\otimes\Lambda(n_j)$.
Then $\fder\fs=\oplus\left((\fder\fs_j)\otimes\Lambda(n_j)\subplus
1\otimes \fvect(n_j)\right)$.  Let $\fg$ be a subalgebra of
$\fder\fs$ containing $\fs$.

1) If the projection of $\fg$ on $1\otimes\fvect(n_j)_{-1}$
coincides with $\fvect(n_j)_{-1}$ for each $j=1, \dots , k$, then
$\fg$ is semi-simple.

2) All semi-simple Lie superalgebras arise in the manner indicated.

Let $\dim\fr=0$, i.e., $\fg$ is semi-simple.  Since $\fg$ is not
almost simple, then, due to Kac's description, the alternative
arises: either $\fg$ contains an ideal $\fii$ of the form
$$
\fii=\fs\otimes \Lambda (n)\quad\text{with simple $\fs$ and
$n>0$}\eqno{(4.1)}
$$
or $\fg=\mathop{\oplus}\limits_{j\leq k}\fs_j$, where each
$\fs_j$ is almost simple and $k>1$.

\sssbegin{Lemma} If $\fg=\mathop{\oplus}\limits_{j\leq k}\fs_j$ and
$k\geq 2$, then any irreducible faithful representation of $\fg$ is
of type B with respect to any its ideal $\fs_j$.\end{Lemma}

\begin{proof} Since the stabilizer of any irreducible representation
of $\fs_j$ is the whole $\fg$, the type of any irreducible
representation of $\fg$ with respect to $\fs_j$ can be either A or B.
Due to faithfulness case A is ruled out.
\end{proof}

\sssbegin{Lemma} Let $\fs$ be a simple Lie superalgebra and
$\fii=\fs\otimes\Lambda(n)$, $n>0$.  Then $\fii$ has no faithful
irreducible finite dimensional representations.
\end{Lemma}

For proof see sec. 4.7.4.

\parbegin{Corollary} If $\fg$ contains an ideal $\fii$ of the form
$(4.1)$, then $\fg$ can not have any faithful irreducible finite
dimensional representation of type A with respect to the ideal
$\fii$.
\end{Corollary}

\parbegin{Corollary} Lemmas  $4.7.1$, $4.7.2$ and Corollary $4.7.2.1$
prove Lemma $4.6$ for semi-simple Lie superalgebras.
\end{Corollary}

Recall the following well-known and simple statement.

\sssbegin{Lemma} If $\fs$ is a simple Lie superalgebra, then
$$
[\fs_{\od}, \fs_{\od}]=\fs_{\ev}\quad {\it and}\quad [\fs_{\ev},
\fs_{\od}]=\fs_{\od}\ .
$$
\end{Lemma}

\ssec{Proof of Lemma 4.7.2} For $n> 0$ the Lie superalgebra $\fii$
contains supercommutative ideal
$$
\fn=\fs\otimes\xi_1\dots\xi_n.
$$
Let us show that $\fn$ is contained in the kernel of any irreducible
representation $\rho$ of $\fii$. Consider a nilpotent ideal
$$
{\fm=\fs\otimes \mathop{\oplus}\limits_{i\ge 1}\Lambda^i(\xi)\subset
\fii.}
$$
As follows from \cite{K}, \cite{Ser1}, any irreducible
finite dimensional representation
of $\fm$ is given by a character $\lambda\in\fm^*$ that vanishes on
$[\fm_{\ev},\fm_{\ev}]\oplus\fm_{\od}$.

For $n=2k$ we have
$$
[\fm_{\ev},\fm_{\ev}]\supset
[\fs_{\od}\otimes\Lambda^1(\xi),\fs_{\od}\otimes\Lambda^{2k-1}(\xi)]
=[\fs_\od,\fs_\od]\otimes \Lambda^{2k}(\xi) =\fn_\ev
$$ and
$\fm_\od\supset\fs_{\od}\otimes\Lambda^{2k}(\xi)=\fn_\od$. For
$n=2k+1>1$ we have
$$
[\fm_{\ev},\fm_{\ev}]\supset
[\fs_{\od}\otimes\Lambda^1(\xi),\fs_{\ev}\otimes\Lambda^{2k}(\xi)]
=[\fs_\od,\fs_\ev]\otimes \Lambda^{2k+1}(\xi) =\fn_\ev $$ and
$\fm_\od\supset\fs_{\ev}\otimes\Lambda^{2k+1}(\xi)=\fn_\od$.  Thus,
$\fn$ is contained in the kernel of any irreducible finite dimensional
representation of the ideal $\fm$,
hence, in the kernel of any
irreducible finite dimensional representation of $\fii$.

If $n=1$ and $\tau$ is an arbitrary irreducible subrepresentation of
$\rho|_{\fn}$, then $\dim \tau=1$ or $\eps$ because $\fn$ is
supercommutative.  Hence, $\tau|_{\fn_\od}=0$.  Heading 1) of
Theorem 4.1.1 implies, that the restriction of $\tau$ onto
$[\fg_\ev\otimes 1,\fn]\supset [\fs_\ev,\fs_\od]\otimes\xi=
\fs_\od\otimes\xi=\fn_\ev$ must vanish.  Thus, $\tau|_{\fn}=0$ and
by heading 3) of Theorem 4.1.1, $\fn\subset\ker \rho$.

Lemma 4.7.2 is proved. \qed

\ssec{Proof of Lemma 4.6 for central extensions of semi-simple Lie
superalgebras} In this section we assume that $\dim \fr=1$, hence
(see Remark 4.3), $\fr$ is the center of $\fg$, and
$\tilde\fg=\fg/\fr$ is semi-simple.

First, consider the case when $\tilde\fg=
\mathop{\oplus}\limits_{i\leq k}\fs_i$, where the $\fs_i$ are almost
simple and $k>1$.  Let $\pi:\fg\tto\tilde\fg$ be the natural
projection.  The Lie superalgebra $\pi^{-1}(\fs_1)=\fii$ is an ideal
in $\fg$ and $\dim \fii>1$.

\sssbegin{Lemma} If $k>1$, then $\rho$ can not be irreducible of
type A with respect to $\fii$.
\end{Lemma}

\begin{proof} Assume the contrary, let $\rho|_{\fii}$ be irreducible.
Set $\fg_+=\mathop{\oplus}\limits_{i> 1}(\fs_i)_{\ev}$ and $\fg_+$
is a Lie algebra.  Then $[\fs_1, \fg_+ ]=0$.  Since $\rho(\fr)$ acts
by scalar operators and $\rho$ is a finite dimensional
representation, it follows that $[\pi^{-1}(\fg_+), \fii]=0$.  Since
$\rho|_{\fii}$ is irreducible by the hypothesis, this means that
$\rho(\pi^{-1} (\fg_+))=\Cee\cdot 1$ and, since $\rho$ is faithful,
this implies $\fg_+=0$.  \end{proof}

\ssec{$\tilde\fg$ contains an ideal $\fii$ of the form $(4.1)$}  The
central extension is defined by a cocycle $c: \tilde\fg\times
\tilde\fg\longrightarrow\Cee$.  The cocycle condition is
$$
c(f, [g, h])=c([f, g], h)+(-1)^{p(f)p(g)}c(g, [f, h])\quad \text{for
any} f, g, h\in\tilde\fg.  \eqno{(4.2)}
$$
As earlier, we assume that $\fg$ has a faithful finite dimensional
representation; so the restriction of $c$ to $\tilde\fg_{\bar
0}\times \tilde\fg_{\bar 0}$ is trivial.  Besides,
$c|_{\tilde\fg_{\bar 0}\times \tilde\fg_{\bar 1}}=0$ by parity
considerations.  Therefore, nonzero values of the cocycle $c$ are
only possible on $\tilde\fg_{\bar 1}\times \tilde\fg_{\bar 1}$.

\sssbegin{Lemma} Let, as above, $\fn=\fs\otimes\Lambda^n(\xi)$. Then
$c|_{\fn\times\tilde\fg}=0$.
\end{Lemma}

\begin{proof} First, let us prove that
$c|_{\fn_{\od}\times\fn_{\od}}=0$.  If $n=2k+1$, it suffices to check
condition $(4.2)$ for the triple
$f=f_{\bar0}\otimes\xi_1\cdot\dots\cdot\xi_n$, $g=g_{\bar
1}\otimes\xi_1\cdot\dots\cdot\xi_n$ and $h=h_{\od}\otimes 1$, where
$f_{\bar 0}\in \fs_{\ev}$ and $g_{\od}, h_{\od}\in \fs_{\od}$.  With
Lemma 4.7.3, the left hand side of $(4.2)$ gives us the values of $c$
on an arbitrary pair of elements from $\fn_{\od}$, whereas the right
hand side vanishes because $[f,g]=0$ and $c(g,[f,h])=0$ since
$p(g)=\ev$.

If $n=2k$, similar arguments are applicable to the triple $f=f_{\bar
1}\otimes\xi_1\cdot\dots\cdot\xi_n$, $g=g_{\bar
0}\otimes\xi_1\cdot\dots\cdot\xi_n$ and $h=h_{\od}\otimes 1$, where
$f_{\bar 1}, h_{\od}\in \fs_{\od}$ and $g_{\ev}, \in \fs_{\ev}$.

Now set $\fL^k=\fs\otimes\Lambda^k(\xi_1, \dots ,\xi_n)$ and let us
verify that $c|_{\fn\times\fL^{k}}=0$.  We will perform the inverse
induction on $k$.  For $k=n$ we have already verified the fact.

Let the statement be true for all $k> k_0$.  Let us show that it is
true for $k=k_0$ as well.  Observe that due to description of
semi-simple Lie superalgebras, $\tilde\fg$ contains $n$ elements
$\eta _i$ such that $\ad\eta
_i|_{\fs\otimes\Lambda(n)}=\partial_{\xi_{i}}+D_i+X_i\otimes\alpha_i$,
where $D_i\in\fvect(0|n)$ and $D_i(0)=0$, $X_i\in\fder\fs,
\alpha_i\in\Lambda(\xi)$.  Let $\varphi\in \Lambda^{k_{0}+1}(\xi_1,
\dots, \xi_n), \psi\in \Lambda^{n}(\xi_1, \dots, \xi_n);\quad g,
h\in\fs$ and $p(g\otimes \varphi)=\bar 0$; $p(h\otimes \psi)=\bar
1$.

Then we have
$$
\renewcommand{\arraystretch}{1.4}
\begin{array}{l}
c([\eta_i, g\otimes\varphi],
h\otimes\psi)=(-1)^{p(g)}\left(c(g\otimes\frac{\partial\varphi}{\partial\xi_{i}},
h\otimes\psi)+c(g\otimes D_{i}\varphi, h\otimes\psi)\right)
+\\
(-1)^{p(g)p(\alpha_i)}c([[X_i,g]\otimes \alpha_i\varphi,
h\otimes\psi).
\end{array}
$$
As $D_i\varphi,\; \alpha_i\varphi\in\mathop{\oplus}\limits_{l\geq
k_{0}+1}\Lambda^l(\xi)$, the last two summands vanish by the inductive
hypothesis.

On the other hand, due to $(4.2)$ we have
$$
c([\eta_i, g\otimes\varphi],h\otimes\psi)=c(\eta_i, [g\otimes\varphi,
h\otimes\psi])+c([\eta_i,  h\otimes\psi], g\otimes\varphi).
$$
As $\deg\varphi+\deg\psi=k_0+1+n>n$, the bracket in the first summand
above vanishes.  The second summand vanishes by parity considerations.
So $c(g\otimes\frac{\partial\varphi}{\partial\xi_{i}},
h\otimes\psi)=0$, for arbitrary
$g\otimes\frac{\partial\varphi}{\partial\xi_{i}}\in\fL_{\bar
1}^{k_{0}}$ and $h\otimes \psi \in\fL_{\bar 1}^{n}=\fn_\od$.
\end{proof}

\sssbegin{Lemma} $\fii$ has no faithful irreducible finite
dimensional representations.
\end{Lemma}

\begin{proof} Word-for-word proof of Lemma 4.7.2 with the
help of Lemma  4.8.2. \end{proof}

\ssbegin{Corollary} $\rho$ can not be of type A with respect to the
ideal $\fii$.\end{Corollary}

\ssec{Summing up} Lemmas 4.2 (Lemma B), 4.4 (Lemma C), 4.5 (Lemma
A), sec. 4.4 and Lemma 4.6 put together prove Main Theorem.  \qed

\bibliographystyle{plain}

\end{document}